\documentclass{article}[10pt]
\usepackage{latexsym}
\usepackage{amsthm}
\usepackage{amssymb}
\newtheorem{theorem}{Theorem}
\newtheorem{definition}{Definition}
\newtheorem{proposition}{Proposition}

\newtheorem{corollary}{Corollary}

\newcommand{\Z}{\mathbb{Z}}

\newcommand{\R}{\mathbb{R}}

\title{A new class of hypercomplex analytic cusp forms}
\author{
D.~Constales \thanks{Department of
Mathematical Analysis, Ghent University, Building
S-22; and Laboratory for Chemical Technology, Ghent University, Building S-5; both at Krijgslaan 281, B-9000 Gent, Belgium. E-mail:
{\tt Denis.Constales@gmail.com}}
\and D. Grob \thanks{Lehrstuhl A f\"ur Mathematik, RWTH Aachen, D-52056 Aachen. E-mail:
{\tt dennisgrob@mathA.rwth-aachen.de}}
\and R. S. Krau{\ss}har \thanks{Arbeitsgruppe Algebra, Geometie und Funktionalanalysis, Fachbereich Mathematik, Technische Universit\"at Darmstadt, Schlo{\ss}gartenstra{\ss}e 7, D-64289 Darmstadt. E-mail:
{\tt krausshar@mathematik.tu-darmstadt.de}}
\and John Ryan \thanks{Department of Mathematics, University of Arkansas, Fayetteville, AR 72701, USA {\tt jryan@uark.edu}}}
 
\begin{document}
\maketitle

\begin{abstract}
In this paper we deal with a new class of Clifford algebra valued automorphic forms on arithmetic subgroups of the Ahlfors-Vahlen group. The forms that we consider are in the kernel of the operator $D \Delta^{k/2}$ for some even $k \in {\mathbb{Z}}$. They will be called $k$-holomorphic Cliffordian automorphic forms. $k$-holomorphic Cliffordian functions are well equipped with many function theoretical tools. Furthermore, the real component functions have also the property that they are solutions to the homogeneous and inhomogeneous Weinstein equation. 
This function class includes the set of $k$-hypermonogenic functions as a special subset. While we have not been able so far to propose a construction for non-vanishing $k$-hypermonogenic cusp forms for $k \neq 0$, we are able to do so within this larger set of functions. After having explained their general relation to hyperbolic harmonic automorphic forms we turn to the construction of Poincar\'e series.  These provide us with non-trivial examples of cusp forms within this function class. Then we establish a decomposition theorem of the spaces of $k$-holomorphic Cliffordian automorphic forms in terms of a direct orthogonal sum of the spaces of $k$-hypermonogenic Eisenstein series and of $k$-holomorphic Cliffordian cusp forms. 
\end{abstract}

\section{Introduction}

The theory of higher dimensional Maa{\ss} wave forms has become a major topic of study in analytic number theory. Maa{\ss} wave forms are automorphic forms that are complex-valued eigensolutions to the Laplace-Beltrami operator 
\begin{equation}
\label{LB}
\Delta_{LB} = x_n^2 (\sum_{i=0}^n\frac{\partial^2 }{\partial x_n^2}) - (n-1) x_n \frac{\partial }{\partial x_n}.
\end{equation}
The classical setting is $n+1$-dimensional upper half-space in the framework of the action of discrete arithmetic subgroups of the orthogonal group. Its study was initiated in 1949 by H. Maa{\ss} in \cite{Maa2}. In the late 1980s this study had a major boost by breakthrough works of J.~Elstrodt, F.~Grunewald and J.~Mennicke, \cite{EGM85,EGM88,EGM90}, A.~Krieg \cite{Kri2,Kri88}, V. Gritsenko \cite{Grit} among many others.  

\par\medskip\par 
 
In the recent book~\cite{KraHabil} another class of automorphic forms on these arithmetic groups has been considered. The context is again $n+1$-dimensional upper half-space. The classes of automorphic forms considered in \cite{KraHabil}, however,  have different analytic and mapping properties. They are null-solutions to the Euclidean Dirac operator $D:=\sum_{i=0}^n \frac{\partial }{\partial x_i} e_i$ (or more generally to higher order iterates of the Dirac operator)  and, in general, they take values in real Clifford algebras.  These are called monogenic (respectively polymonogenic) automorphic forms.  The monogenic automorphic forms on upper half-space in turn can be embedded into the general framework of $k$-hypermonogenic automorphic forms. 
This has been done in our previous work \cite{CKR2007}. The class of $k$-hypermonogenic functions, discussed for example in  \cite{Leutwiler,HDiss,slberlin,slnewnew} contains the set of monogenic functions ($k=0$). Monogenic functions were studied by numerous authors. For an overview about their associated function theory see for instance \cite{DSS,GS} and elsewhere. $k$-hypermonogenic functions are solutions to a modification of the Dirac operator associated with the hyperbolic metric on upper half-space. 

\par\medskip\par

As explained in \cite{CKR2007}, the more general framework of $k$-hypermonogenic automorphic forms allows us to link Maa{\ss} wave forms with the class of monogenic automorphic forms.  In particular, it turned out that one can construct Maa{\ss} wave forms from monogenic automorphic forms. 

However, for $k \neq 0$ we have only been able to construct $k$-hypermonogenic Eisenstein series so far. We were not able to propose a construction for non-vanishing $k$-hypermonogenic cusp forms. The idea to extend the classical construction 
$$
P(z,w) = \sum_{M  \in SL(2,\mathbb{Z})} (cz+d)^{-k}\Bigg(w+\frac{az+b}{cz+d}\Bigg)^{-k},\; Im(z) > 0,\;Im(w) > 0
$$
(where $M= \left( \begin{array}{rr} a & b \\ c & d \end{array} \right)$) fails in this context. The reason is that a translation induced by $w$ into the upper half-space $x_n$ direction, does not preserve $k$-hypermonogenicity. 

This motivates us to consider a larger class of functions that contains the class of $k$-hypermonogenic functions as a special subset, but which is still sufficiently specialized so that one still has many function theoretical tools available. In the cases where $k \in {\mathbb{Z}}$ and $k$ even, a natural extension of $k$-hypermonogenic functions in that sense are functions that are annihilated by the operator $D \Delta^{k/2}$. Every $k$-hypermonogenic function is an element of Ker $D \Delta^{k/2}$. However, the converse is not true. For negative $k$, the Laplacian is formally substituted by the product of the Teodorescu operator with its conjugate. The Teodorescu operator is the right inverse to the Dirac operator $D$. However, we can directly relate the operator $D \Delta^{-k/2}$ with the differential operator $D\Delta^{(k+2)/2}$. By multiplying an element from Ker $D \Delta^{-k/2}$ with the simple scaling factor $\frac{1}{x_n^k}$, we again get an element from the kernel of $D\Delta^{(k+2)/2}$. Therefore, it suffices to consider either non-negative or non-positive even numbers $k$.   
In the special case where $k=n-1$ and where $n$ is odd we deal with the holomorphic Cliffordian functions that were considered earlier in \cite{LR,LR2,LR3,Ra2008} and elsewhere. As shown for instance in \cite{LR}, the set of holomorphic Cliffordian functions actually is endowed with many function theoretical tools that are also offered for complex holomorphic functions. This includes for instance a Cauchy type integral formula, a residue calculus, Taylor and Laurent expansions and much more. For general $k$ we naturally call the functions from Ker $D \Delta^{k/2}$ $k$-holomorphic Cliffordian functions. In fact, those function theoretical tools turned out to be very useful in the solution of boundary value problems from harmonic analysis on the quotient manifold or orbifold that arise by factoring upper half-space modulo  arithmetic subgroups of the orthogonal group. See for instance \cite{paper4}, in which we used monogenic automorphic forms to solve boundary value problems related to the Laplacian on these quotient structures.  

In turn the real components of a $k$-holomorphic Cliffordian function satisfy the homogeneous Weinstein equation or the inhomogeneous Weinstein equation of parameter $-k$, respectively. This is true as well for positive $k$ but also for negative ones. The Weinstein equation is a second order partial differential equation.  
   
\par\medskip\par

In Section~2 we recall some basic notations on Clifford algebras and associated arithmetic subgroups. In Section~3 we introduce differential operators that generalize the classical Cauchy-Riemann operator to higher dimensional Euclidean and hyperbolic spaces. Then we explain how $k$-holomorphic Cliffordian functions are related to the Weinstein equation and how the solutions to the Weinstein equation can be characterized in terms of nullsolutions to the pluriharmonic equation $\Delta^{(k+2)/2}f =0$.  

\par\medskip\par

In Section~4 we introduce $k$-holomorphic Cliffordian automorphic forms and discuss some basic properties of them. We show how we can construct $k$-holomorphic Cliffordian Eisenstein series from $k$-hypermonogenic Eisenstein series. Then we carefully discuss the relation between $k$-holomorphic Cliffordian automorphic forms and Maa{\ss} wave forms. We explain how the examples of Eisenstein- and Poincar\'e series from \cite{EGM88,EGM90,Kri2,Kri88} are related to the class of $k$-holomorphic Cliffordian automorphic forms. 

\par\medskip\par

In Section~5, which is the main section of this paper, we give an explicit construction of $k$-holomorphic Cliffordian Poincar\'e series. These provide us with non-trivial examples of cusp forms for infinitely many $k$. These Poincar\'e series are not $k$-hypermonogenic. 
 
Then we compute the Fourier expansion of $k$-holomorphic Cliffordian automorphic forms. With this tool in hand, we finally establish an important decomposition of the space of Clifford algebra valued $k$-holomorphic Cliffordian automorphic forms. This decomposition is expressed in terms of a direct orthogonal sum of the space of $k$-hypermonogenic Eisenstein series and the space of $k$-holomorphic Cliffordian cusp forms. 

The orthogonality is established with respect to the same Petersson inner product considered in the context of Maa{\ss} wave forms in $n+1$ real dimensions and is a nice analogy of the orthogonality result proved in \cite{EGM90} for Maa{\ss} wave forms.   

\par\medskip\par

These results also provide analogues to the classical theory
of complex analytic automorphic forms. See for example \cite{Freitag,schoeneberg}. Furthermore, they give some first insight into the fundamental structure of the space of $k$-holomorphic Cliffordian automorphic forms. It is an amazing fact that the space of $k$-holomorphic Cliffordian cusp forms (none of them is $k$-hypermonogenic) provides us exactly with the orthogonal complement of the space of $k$-hypermonogenic Eisenstein series. 

\section{Preliminaries}

\subsection{Clifford algebras}

We introduce the basic notions of real Clifford algebras over the Euclidean space $\mathbb{R}^n$. For details, see for instance~\cite{DSS}. Throughout this paper, $\{e_1,\ldots,e_n\}$ stands for the standard orthonormal basis in the Euclidean space $\mathbb{R}^n$ and $Cl_n$ denotes its associated real Clifford algebra
in which the relation $e_i e_j + e_j e_i = - 2 \delta_{ij}$ holds. Each vector $\vec{x} \in \mathbb{R}^n \backslash\{0\}$ has an inverse element, given by $-\vec{x}/\|\vec{x}\|^2$. A basis for
the algebra $Cl_{n}$ is given by $1, e_{1},\ldots,e_{n},\ldots,e_{j_{1}}\ldots e_{j_{r}},$ $\ldots, e_{1}\ldots e_{n}$ 
where $j_{1}<\ldots j_{r}$ and $1\leq r\leq n$.
The reversion anti-automorphism is defined by  $\widetilde{ab} = \tilde{b} \tilde{a}$, where $\widetilde{e_{j_{1}}\ldots e_{j_{r}}} = e_{j_{r}}\ldots e_{j_{1}}$. The conjugation anti-automorphism is defined by $\overline{ab}=\overline{b} \;\overline{a}$ where $\overline{e_i}=-e_i$ for $i=1,\ldots, n$. 
The main involution is defined by $(ab)^{'} = a^{'}b^{'}$ and $e_i^{'}=-1$ for $i=1,\ldots,n$.
Furthermore, we also need the following automorphism $*: Cl_n \rightarrow Cl_n$ defined by the
relations: $e_n^{*}=-e_n$, $e_i^{*}= e_i$ for $i=0,1,\ldots,n-1$ and $(ab)^{*}=a^{*}b^{*}$. Any element $a \in Cl_n$ may be uniquely decomposed in the form $a = b + c e_n$, where $b$ and $c$ belong to $Cl_{n-1}$. Based on this definition one defines the mappings $P: Cl_n \rightarrow Cl_{n-1}$ and $Q: Cl_n \rightarrow Cl_{n-1}$ by $Pa=b$ and $Qa=c$.

An important subspace of the Clifford algebra is the space of paravectors $\mathbb{R} \oplus \mathbb{R}^n$ which consists of elements of the form $x=x_0+x_1e_1 + \cdots + x_n e_n=x_0+\vec{x}$. Every non-zero paravector has an inverse of the form $x^{-1}=\frac{\overline{x}}{|x|^2}$, where $|\cdot|$ stands for the Euclidean norm on $\mathbb{R} \oplus \mathbb{R}^n$. On the whole Clifford algebra one can introduce a pseudo norm by $\|a\| = ( \sum\limits_A |a_A|^2)^{1/2}$, where $a_A$ are the scalar components of a Clifford number $a=a_0 + a_1 e_1 + \cdots + a_n e_n + a_{12} e_1 e_2 + \cdots + a_{1...n} e_1 \cdots e_n$. For simplicity we represent the Clifford number $a$ in the form $a = \sum\limits_{A \subseteq P(\{1,...,n\})} a_A e_A$ where $P(\{1,...,n\})$ is the potential set of the set $\{1,\ldots,n\}$ and $e_A$ represents the product of all elements $e_i$ where $i \in A$.   

\subsection{Discrete arithmetic subgroups of $GAV(\mathbb{R} \oplus\mathbb{R}^{n})$}

In this paper we deal with Clifford algebra-valued functions that are quasi-invariant under arithmetic subgroups of the general Ahlfors-Vahlen group that act totally discontinuously on upper half-space $H^{+}(\mathbb{R} \oplus\mathbb{R}^n) =\{x \in \mathbb{R} \oplus \mathbb{R}^n\;:\;x_n > 0\}$. These can be regarded as generalizations of the classical holomorphic modular forms in the context of monogenic, $k$-hypermonogenic and $k$-holomorphic Cliffordian functions. These particular function classes will be described in the following section, after having outlined the algebraical setting.  

\par\medskip\par

Arithmetic subgroups of the special Ahlfors-Vahlen group that act totally discontinuously on upper half-space were for instance considered in \cite{Maa2,EGM87,EGM90}. 

The set that consists of Clifford valued matrices $\left(\begin{array}{cc} a & b \\ c & d \end{array}\right)$ whose coefficients satisfy the conditions below forms a group under matrix multiplication. It is called the general Ahlfors-Vahlen group, $GAV(\mathbb{R} \oplus \mathbb{R}^n)$. The action of the general Ahlfors-Vahlen group on $\mathbb{R} \oplus \mathbb{R}^n$ is described by the associated M\"obius transformation.  
Following for example \cite{a,EGM87}, M\"obius transformations in $\mathbb{R} \oplus\mathbb{R}^n$ can be
represented as
$$
T: \mathbb{R} \oplus\mathbb{R}^n \cup \{\infty\} \rightarrow \mathbb{R} \oplus\mathbb{R}^n \cup \{\infty\},\; T(x)=(ax+b)(cx+d)^{-1}
$$
with coefficients $a,b,c,d$ from $Cl_n$ satisfying
\begin{enumerate}
\item[(i)] \ $a,b,c,d$ are products of paravectors  
\item[(ii)] \ $a\tilde{d} - b\tilde{c} \in \R \setminus \{ 0 \} $ 
\item[(iii)] \ $ac^{-1}, \, c^{-1}d \in \R^{n+1} $ \, for \, $ c \neq 0 \, $ and  \, $ bd^{-1} \in \R^{n+1} $ \, for \,$ c=0 $ .
\end{enumerate}
We also describe the group action in the way $GAV(\mathbb{R} \oplus \mathbb{R}^n) \times H^+(\mathbb{R} \oplus \mathbb{R}^n) \to H^+(\mathbb{R} \oplus \mathbb{R}^n)$ by  $(M,x)\mapsto M<x> = (ax+b)(cx+d)^{-1}$. 

\par\medskip\par

The subgroup consisting of those matrices from $GAV(\mathbb{R} \oplus\mathbb{R}^n)$ that  satisfy $ad^*-bc^*=1$ is called the special Ahlfors-Vahlen group. It is denoted by $SAV(\mathbb{R} \oplus\mathbb{R}^n)$. The automorphism group of upper half-space $H^+(\mathbb{R} \oplus \mathbb{R}^n)$ is the group $SAV(\mathbb{R} \oplus\mathbb{R}^{n-1})$. 

\par\medskip\par

Next we recall the definition of the rational Ahlfors-Vahlen group acting on $H^+(\mathbb{R} \oplus\mathbb{R}^n)$.
\begin{definition}
The rational Ahlfors-Vahlen group $SAV(\mathbb{R} \oplus\mathbb{R}^{n-1},\mathbb{Q})$ is the set of matrices
$\left(\begin{array}{cc} a & b \\ c & d \end{array}\right)$ from $SAV(\mathbb{R} \oplus\mathbb{R}^{n-1})$ that satisfy \\[0.2cm]
$ (i)\quad a \overline{a}, b \overline{b}, c \overline{c}, d \overline{d} \in {\mathbb{Q}},\\
  (ii) \quad a\overline{c},b\overline{d} \in {\mathbb{Q}}
\oplus \mathbb{Q}^{n},\\
  (iii)  \quad a x \overline{b} + b \overline{x} \;\overline{a}, c x \overline{d} + d \overline{x} \;\overline{c} \in    {\mathbb{Q}}\quad (\forall x \in \mathbb{Q} \oplus 
\mathbb{Q}^{n})\\
 (iv) \quad a x \overline{d}+b \overline{x}\;\overline{c}  
\in {\mathbb{Q}}
\oplus \mathbb{Q}^{n}
\quad (\forall x \in \mathbb{Q} \oplus \mathbb{Q}^{n})$.
\end{definition}
Next we need
\begin{definition}
A $\mathbb{Z}$-order in a rational Clifford algebra is a subring $R$ such that the additive group of $R$ is finitely generated and contains a $\mathbb{Q}$-basis of the
Clifford algebra.
\end{definition}
The following definition, cf. \cite{EGM90}, provides us with a whole class of arithmetic subgroups of the Ahlfors-Vahlen group which act totally discontinuously on upper half-space. 
\begin{definition}
Let ${\cal{I}}$ be a ${\mathbb{Z}}$-order in ${Cl}_{n}$ which is stable under the reversion and the main involution ${'}$ of ${Cl}_{n}$. Then
$$
\Gamma_{n-1}({\cal{I}}) := SAV(\mathbb{R} \oplus \mathbb{R}^{n-1},{\mathbb{Q}}) \cap Mat(2,{\cal{I}}).
$$
For an $N \in {\mathbb{N}}$ the principal congruence subgroup of
$SAV(\mathbb{R} \oplus\mathbb{R}^{n-1},{\cal{I}})$ of level $N$ is defined by
$$
\Gamma_{n-1}({\cal{I}})[N] := \Big\{ \left(\begin{array}{cc} a & b\\ c & d \end{array}\right) \in \Gamma_{n-1}({\cal{I}}) \;\Big|\; a-1,b,c,d-1 \in N {\cal{I}}\Big\}.
$$
\end{definition}
Notice that all the groups $\Gamma_{n-1}({\cal{I}})[N]$ have finite index in $\Gamma_{n-1}({\cal{I}})$. Therefore, all of them are discrete groups and act totally discontinuously on upper half-space $H^{+}(\mathbb{R} \oplus{\mathbb{R}}^n)$. 
The proof of the total discontinuous action can be done in the same way as in \cite{Freitag}. 

\par\medskip\par

{\bf Special notation}. 
In the sequel let us denote the fundamental domain of $\Gamma_{n-1}({\cal{I}})[N]$ by ${\cal{F}}_{n-1}[N]$. 
Furthermore, denote the subgroup of translation matrices contained in $\Gamma_{n-1}({\cal{I}})[N]$ by ${\cal{T}}_{n-1}({\cal{I}})[N]$. The associated $n$-dimensional lattice in $\mathbb{R} \oplus \mathbb{R}^{n-1}$ will be denoted by $\Lambda_{n-1}({\cal{I}})[N]$ and its dual by $\Lambda_{n-1}^{*}({\cal{I}})[N]$. The dual lattice is contained in the subspace $\mathbb{R} \oplus \mathbb{R}^{n-1}$, too. The standard fundamental period cell of $\Lambda_{n-1}({\cal{I}})[N]\subset \mathbb{R} \oplus\mathbb{R}^{n-1}$ will be denoted by $L_{n-1}[N]$. 

\par\medskip\par

The simplest concrete examples for $\Gamma_{n-1}({\cal{I}})$ are obtained by taking for ${\cal{I}}$ the standard ${\mathbb{Z}}$-order in the Clifford algebras $Cl_p$, i.e.,  ${\cal{O}}_p := \sum_{A \subseteq  P(1,\ldots,p)} \mathbb{Z} e_A$ where $p \le  n-1$.  
In this case, the group $\Gamma_{n-1}({\cal{I}})$ coincides with the special hypercomplex modular group of parameter $n-1$ considered in \cite{KraHabil}. The special hypercomplex modular group of parameter $p$ with $p \in \{0,... n-1\}$ is generated by the matrices 
$$
J:= \left(\begin{array}{cc} 0 & -1 \\ 1 & 0 \end{array}\right),T_{1} := \left(\begin{array}{cc} 1 & 1 \\ 0 & 1  \end{array}\right), T_{e_1} := \left(\begin{array}{cc} 1 & e_1 \\ 0 & 1  \end{array}\right),\ldots, T_{e_p}:= \left(\begin{array}{cc} 1 & e_p \\ 0 & 1  \end{array}\right).
$$
In this particular case, we have ${\cal{T}}_{p}({\cal{I}}) = <T_1,T_{e_1},\ldots,T_{e_p}>$. The associated period lattice is the orthonormal $p+1$-dimensional lattice $\mathbb{Z} + \mathbb{Z} e_1 + \cdots + \mathbb{Z} e_p$ which is self-dual. Its standard fundamental period cell is $[0,1]^{p+1}$. 

\par\medskip\par

In the case $n=3$ the upper half-space $\mathbb{R} \oplus \mathbb{R}^3$  can be identified with the upper quaternionic half-space. In this setting further important examples for ${\cal{I}}$ are the quaternionic orders, in particular the Hurwitz order, considered for example in \cite{Kri2,Kri88,Reid}.  

\section{Differential operators}

\subsection{Dirac type operators over $n+1$-dimensional Euclidean and hyperbolic spaces}

{\bf Monogenic functions.} Let $U \subseteq \mathbb{R} \oplus\mathbb{R}^n$ be an open set. Then a real differentiable function $f: U \rightarrow Cl_n$ that  satisfies $Df = 0$, respectively $fD = 0$, where $$D := \frac{\partial }{\partial x_0}  + \frac{\partial }{\partial x_1} e_1 + \cdots + \frac{\partial }{\partial x_n}e_n$$ is the Euclidean Dirac operator, is called left monogenic (respectively right monogenic), cf. for instance \cite{DSS}. Due to the non-commutativity of $Cl_{n}$ for $n>1$, both
classes of functions do not coincide with each other. However, $f$ is left monogenic if and only if $\tilde{f}$ is right monogenic. The Dirac operator factorizes the Euclidean Laplacian $\Delta = \sum_{j=0}^n \frac{\partial^2}{\partial x_j^2}$, viz $\overline{D} D = \Delta$. Every real component of a monogenic function hence is harmonic.  

\par\medskip\par

An important property of the $D$-operator is is quasi-invariance under M\"obius transformations acting on the complete Euclidean space $\mathbb{R} \oplus\mathbb{R}^n$.

\par\medskip\par

Let $M \in GAV(\mathbb{R} \oplus\mathbb{R}^n)$ and $f$ be a left monogenic function in the variable
$y=M<x>=(ax+b)(cx+d)^{-1}$. Then, following, e.g. \cite{r85}, the function $$g(x):=\frac{\overline{cx+d}}{\|cx+d\|^{n+1}}f(M<x>)$$ is again left  monogenic in the variable $x$ for any $M \in GAV(\mathbb{R} \oplus\mathbb{R}^n)$.

\par\medskip\par

{\bf $k$-hypermonogenic functions.}
The class of monogenic functions belongs to the more general class of so-called $k$-hypermonogenic functions. These are defined as the null-solutions to the system
$$
D f + k \frac{(Q f)^{'}}{x_n} = 0
$$
where $k \in \mathbb{R}$. For the basic properties of $k$-hypermonogenic functions, see for example \cite{slberlin,slnew,qber}. 

\par\medskip\par

In the case $k=0$, we are dealing with the set of left monogenic functions introduced earlier. The particular solutions associated to the case $k=n-1$ coincide with the null-solutions to the hyperbolic Hodge-Dirac operator with respect to the hyperbolic metric on upper half space. These are often called hyperbolic monogenic functions or simply hypermonogenic functions, see \cite{Leutwiler}.

\par\medskip\par

Any M\"obius transformation $T(x)=M<x>$ induced by  matrices $M$ belonging to $SAV(\mathbb{R} \oplus\mathbb{R}^{n-1})$ (as well as any of its subgroups) preserves upper half-space. Let $y=M<x>=(ax+b)(cx+d)^{-1}$ be the image of a point $x$ under such a M\"obius transformation $M$. 

If $f$ is $k$-hypermonogenic in the variable $y=M<x>=(ax+b)(cx+d)^{-1}$ then such a solution $f(y)$ is transformed to the $k$-hypermonogenic function by 
\begin{equation}
\label{kinv}
F(x):=\frac{\overline{cx+d}}{\|cx+d\|^{n+1-k}}f(M<x>).
\end{equation}
The new function $F(x)$ satisfies the $k$-hypermonogenic equation in the variable $x$.
See for instance \cite{slnewnew}.

Therefore, in view of~(\ref{kinv}), if $f$ is a function that is left $k$-hypermonogenic on the whole half-space then so is $\frac{\overline{cx+d}}{\|cx+d\|^{n+1-k}} f(M<x>)$.

\par\medskip\par

However, this invariance only holds for matrices from $SAV(\mathbb{R} \oplus\mathbb{R}^{n-1})$ and not for all matrices from $SAV(\mathbb{R} \oplus\mathbb{R}^{n})$. A translation in the argument of a $k$-hypermonogenic function in the $e_n$-direction in general does not result into a $k$-hypermonogenic function again. Consider next only even integer values for $k$. Then the function class of $k$-holomorphic Cliffordian functions that we discuss in the follow up section is an extension of $k$-hypermonogenic functions. $k$-holomorphic Cliffordian functions have the extra property of being invariant under the whole group $SAV(\mathbb{R} \oplus\mathbb{R}^{n})$. A further problem that occurs in the context of $k$-hypermonogenic functions is the following.  Multiplication of $e_n$ from the right to a $k$-hypermonogenic function does in general not give again a $k$-hypermonogenic function. However, the larger function class that we are going to introduce now, has the property that right multiplication of any number from $Cl_n$ with a $k$-holomorphic Cliffordian function remains $k$-holomorphic Cliffordian, see for instance \cite{Ra2008}. These extra properties have crucial advantages. 

\subsection{$k$-holomorphic Cliffordian functions and the Weinstein equation}
To proceed in this direction we introduce 
\begin{definition} Let $n \in \mathbb{N}$ and suppose that $k$ is an even  positive integer. Let $U \subset \mathbb{R} \oplus\mathbb{R}^{n}$ be an open subset. Then we call a function $f:U \to Cl_n$ $k$-holomorphic Cliffordian if   
\[
  D \bigtriangleup^{\frac{k}{2}} f = 0.
\]
\end{definition}
In the particular case $k=n-1$ ($n$ odd) we deal with the class of holomorphic Cliffordian functions, discussed for instance in \cite{LR,LR2,LR3}. As explained in \cite{LR}, simple examples of holomorphic Cliffordian functions are the positive and negative powers of the hypercomplex variable $x=x_0+e_1 x_1 + \cdots +e_n x_n$. One has Taylor and Laurent series expansions in terms of permutational products of powers of the standard hypercomplex variables. Moreover one has close analogues of Cauchy's integral formula and a residue calculus for this function class. In fact many of these tools such as the Cauchy integral formula carry over to the even more general framework of $k$-holomorphic Cliffordian functions. The reason for that is that this function class is a very natural extension of the class of $k$-hypermonogenic functions. For the latter one has a very well established function theory. Before we shed more light on that connection, we will also need to introduce $k$-holomorphic Cliffordian functions for negative even integers $k$. To introduce negative powers of the Laplacian on upper half-space, we make use of the well-known Teodorescu operator, used for instance in \cite{GS} and elsewhere. Adapting from \cite{GKRS}, for the half-space setting the latter is defined by
$$
[Tf](x):= -\frac{1}{A_{n+1}} \int_{H^+(\mathbb{R}\oplus\mathbb{R}^{n})}\Bigg(\frac{\overline{y-x}}{\|y-x\|^{n+1}} - \frac{\overline{y-z}}{\|y-z\|^{n+1}}\Bigg)f(y) dV(y)
$$
where $x$ is an element of the upper half-space and $z$ and element of the lower half-space. Here, $A_{n+1}$ stands for the surface of the unit sphere in an $n+1$-dimensional space. As shown in \cite{GS,GKRS} and elsewhere, this operator is the right inverse of the $D$-operator, i.e. $DT f = f$. In view of the identity $\overline{D}D=\Delta$ we can formally express $\Delta^{-1}$ as $\overline{T}{T}$ on upper half-space. 

\par\medskip\par

This interpretation permits us to also introduce $k$-holomorphic Cliffordian functions for negative even integers $k$. However, as we shall see later, there is a simpler way to describe elements in Ker $D \Delta^{k/2}$ for negative $k$. We will be able to prove that if $f$ is a $k$-holomorphic function, then the function $g:=\frac{f}{x_n^k}$ will be $-k+2$-holomorphic Cliffordian. Since $f$ and $g$ only  differ by the scaling factor $\frac{1}{x_n^k}$, we can treat a $-k$-holomorphic Cliffordian function completely by its associated $k+2$-holomorphic Cliffordian function just differing from it by that scalar factor. So without loss of generality one can restrict oneself to consider non-negative values for $k$.

\par\medskip\par

First of all we start by proving that for any even $k \in \mathbb{Z}$ every $k$-hypermonogenic function is also $k$-holomorphic Cliffordian. In the case $k=n-1$ where $n$ is an odd positive integer, a proof of this statement is sketched in \cite{Leutwiler}, see Theorem~2.5. We can say even more:

\begin{theorem}\label{satz_eins}

(i) Let $k \in \Z$ be even and $f$ a $k$-hypermonogenic function. Then
\[
  D \bigtriangleup^{\frac{k}{2}} f = 0.
\]
(ii) Furthermore, a function $g$ satisfies 
\[
  \bigtriangleup^{\frac{k+2}{2}} g = 0
\]
if and only if $Pg$ satisfies the homogeneous Weinstein equation 
$$
\Bigg(\Delta - \frac{k}{x_n}\frac{\partial }{\partial x_n}\Bigg) (Pg) = 0
$$
and $Qg$ the inhomogeneous Weinstein equation with parameter $-k$, i.e.
$$
\Bigg(\Delta - \frac{k}{x_n}\frac{\partial }{\partial x_n}\Bigg) (Qg) = -\frac{k}{x_n} (Qg).
$$
\end{theorem}
{\bf Remark}: Notice that Ker $D \Delta^{k/2} \subset$ Ker $\Delta^{(k+2)/2}$. 
\begin{proof}
(i) Let $f$ be a solution of
\begin{equation}\label{eins}
  x_n \bigtriangleup f + c \frac{\partial f}{\partial x_n} = c \frac{f}{x_n},
\end{equation}
taken as an equivalent substitute for
\begin{equation}\label{zwei}
   x_n^2 \bigtriangleup f + \frac{c}{x_n} \frac{\partial f}{\partial x_n} = c \frac{f}{x_n^2}.
\end{equation}
We can consider (\ref{eins}) instead of (\ref{zwei}), as we only regard the upper half-space with $x_n > 0$. 
We define a function $v$ by
\[
 v = x_n \bigtriangleup f .
\]
Then we obtain
\begin{eqnarray*}
 \bigtriangleup v &=& \bigtriangleup \left( -c \frac{\partial f}{\partial x_n} + c \frac{f}{x_n} \right) \\
&=& -c \frac{\partial \bigtriangleup f}{\partial x_n} + (\bigtriangleup \frac{c}{x_n})f + 2 \sum_{i=0}^n \frac{\partial}{\partial x_i} \frac{c}{x_n} \frac{\partial}{\partial x_i} f + c \frac{\bigtriangleup f}{x_n} \\
&=& -c \frac{\partial \bigtriangleup f}{\partial x_n} + (2 \frac{c}{x_n^3})f + 2 \frac{-c}{x_n^2} \frac{\partial}{\partial x_n} f + c \frac{\bigtriangleup f}{x_n} \\
&=& -c \frac{\partial \bigtriangleup f}{\partial x_n} + 2 \frac{1}{x_n} \bigtriangleup f + c \frac{\bigtriangleup f}{x_n} \\
&=& -c \frac{\partial \bigtriangleup f}{\partial x_n} + (c+2) \frac{\bigtriangleup f}{x_n}.
\end{eqnarray*}
Further, we get
\begin{eqnarray*}
 \bigtriangleup v &=& \bigtriangleup (x_n \bigtriangleup f) \\
&=& 2 \frac{\partial \bigtriangleup f}{\partial x_n} + \bigtriangleup^2 f. 
\end{eqnarray*} From this equality we can deduce
\[
  2 \frac{\partial \bigtriangleup f}{\partial x_n} + \bigtriangleup^2 f = -c \frac{\partial \bigtriangleup f}{\partial x_n} + (c+2) \frac{f}{x_n},
\]
and thus
\[
   \bigtriangleup^2 f = -(c+2) \frac{\partial \bigtriangleup f}{\partial x_n} + (c+2) \frac{f}{x_n},
\]
so for any solution $f$ of (\ref{eins}) it holds that $w:=\bigtriangleup f$ is a solution of (\ref{eins}), if we replace $c$ by $c+2$. 
As (\ref{eins}) and (\ref{zwei}) are equivalent for $x_n \neq 0$, the proposition remains true for (\ref{zwei}) due to the identity theorem for hypermonogenic functions.\\
Combining this result with Lemma 2.6 from the article \cite{Leutwiler} leads to the following statement.
 
If $f$ is $k$-hypermonogenic for a positive $k$, then $\bigtriangleup f$ is $(k-2)$-hypermonogenic. 
So if $f$ is $k$-hypermonogenic, then $\bigtriangleup^{\frac{k}{2}} f$ is monogenic, and thus $D \bigtriangleup^{\frac{k}{2}} f = 0$.
\par\medskip\par
Now let $T$ be the Teodorescu operator and set $v:=\overline{T}Tf$. Suppose that $f$ is a solution of (\ref{zwei}).
Using the calculations above, we obtain 
\begin{eqnarray*}
 && \bigtriangleup \left( -c \frac{\partial }{\partial x_n} v + c \frac{v}{x_n} \right) \\
&=& -c \frac{\partial }{\partial x_n} \bigtriangleup v + (c+2) \frac{\bigtriangleup v}{x_n} \\
&=& -c \frac{\partial }{\partial x_n} f + (c+2) \frac{f}{x_n}. 
\end{eqnarray*}
Thus, we have
\begin{eqnarray*}
 x_n \bigtriangleup f &=& - c \frac{\partial f}{\partial x_n} + c \frac{f}{x_n} \\
&=& \bigtriangleup \left( -(c-2) \frac{\partial }{\partial x_n} v + (c-2) \frac{v}{x_n} \right) + 2 \frac{\partial f}{\partial x_n} 
\end{eqnarray*}
Using the calculations above, we obtain that 
\begin{eqnarray*}
x_n \bigtriangleup f &=& \overline{T}T \left( \bigtriangleup \left( x_n \bigtriangleup f \right)  \right) \\
&=& \overline{T}T \left( 2 \frac{\partial }{\partial x_n} \bigtriangleup f + \bigtriangleup^2 f \right) \\
&=& 2 \frac{\partial }{\partial x_n} f + \bigtriangleup f
\end{eqnarray*}
Combining both equations with each other leads to
\[
 2 \frac{\partial }{\partial x_n} f + \bigtriangleup f = \bigtriangleup \left( -(c-2) \frac{\partial }{\partial x_n} v + (c-2) \frac{v}{x_n} \right) + 2 \frac{\partial f}{\partial x_n} 
\]
This in turn leads to
\[
  \bigtriangleup^2 v = \bigtriangleup f = \bigtriangleup \left( -(c-2) \frac{\partial }{\partial x_n} v + (c-2) \frac{v}{x_n} \right),
\]
and thus
\[
 \bigtriangleup v =  -(c-2) \frac{\partial }{\partial x_n} v + (c-2) \frac{v}{x_n}.
\]
So, if $f$ is $(-k)$-hypermonogenic, then $\overline{T}T f$ is $(-k+2)$-hypermonogenic and $\overline{T}T^{\frac{k}{2}} f = \bigtriangleup^{-\frac{k}{2}} f$ is monogenic for $k$ being integer valued, even and positive. We have $ D\bigtriangleup^{-\frac{k}{2}} f = 0 $.
\par\medskip\par
(ii) Let us prove the equivalence now:

Suppose that a function $f$ satisfies the homogeneous Weinstein equation in its $P$-part and the inhomogeneous Weinstein equation in its $Q$-part for an even integer $k$. Then we can apply the same calculations as above to show that applying $\bigtriangleup$ or  $\bigtriangleup^{-1}$ makes $k$ increase or decrease, respectively. Thus, we again arrive at the fact that $\bigtriangleup^{\frac{k}{2}} f$ is harmonic and $\bigtriangleup^{\frac{k+2}{2}} f = 0$ \\*[0.2cm]
We have proved that if $f$ is solution of a Weinstein equation with constant $c$, then $\bigtriangleup f$ is a solution for the case $c+2$ and $\bigtriangleup^{-1} f$ for the case $c-2$ (this holds for both the homogeneous and the inhomogeneous Weinstein equation).\\
Now, if $f$ is a null-solution to $\bigtriangleup^{\frac{k+2}{2}}$, then $\bigtriangleup^{\frac{k}{2}} f$ is harmonic and thus $\bigtriangleup^{\frac{k}{2}} f$ is a solution of the Weinstein equation with $c=0$. If $k$ is positive, then $ \bigtriangleup^{-\frac{k}{2}} \bigtriangleup^{\frac{k}{2}} f$ is a solution in the case $c= 2 \cdot \frac{-k}{2}=-k$. The case where $k$ is negative can be treated analogously by applying $\bigtriangleup$. This completes the proof of this theorem.
\end{proof}
As a consequence the $k$-hypermonogenic kernel functions $G_k(x) := \frac{\overline{x}}{\|x\|^{n+1-k}}$ are for all $k \in 2 \mathbb{Z}$ also $k$-holomorphic Cliffordian. 

Notice that not every $k$-holomorphic Cliffordian function is $k$-hypermonogenic. Take for instance $k=n-1$ and suppose that $n$ is odd. Then the functions $x e_n$ and $x + e_n$ (with $x=x_0+e_1 x_1 + \cdots +e_n x_n$) are both $n-1$-holomorphic Cliffordian. However, they are not $n-1$-hypermonogenic as one can verify by a simple calculation.  

\par\medskip\par
Nevertheless many nice properties for $k$-hypermonogenic functions still carry over to the larger class of $k$-holomorphic Cliffordian functions. 

A very important property is that the transformation $f \mapsto \frac{f e_n}{x_n^k}$ which transforms a $k$-hypermonogenic function  into a $-k$-hypermonogenic function, does also transform a $k$-holomorphic Cliffordian function into a $-k+2$-holomorphic Cliffordian function. In fact, one possibly can even establish a stronger result, namely that $f \mapsto \frac{f e_n}{x_n^k}$ is even $-k$-holomorphic Cliffordian. For our needs we will only prove the weaker version of this statement. 

Actually, when we deal with $k$-holomorphic Cliffordian functions the factor $e_n$ on the right-hand side may be omitted, because right multiplication with $e_n$ still preserves $k$-holomorphicity  (but not $k$-hypermonogenicity).

We first prove that 
\begin{theorem}\label{satz_zwei}
 Let $k$ be an even integer. Suppose that $f$ is a sufficiently many times real differentiable function with
\[
   \bigtriangleup^{\frac{k}{2}+1} f = 0.
\]
Then the function $g:=\frac{f}{x_n^k}$ satisfies 
\[
   \bigtriangleup^{\frac{-k}{2}+1} g = 0.
\]
\end{theorem}

\begin{proof}
Since $f$ is a nullsolution to $D \bigtriangleup^{\frac{k}{2}} f = 0$, we have by Theorem~\ref{satz_eins}(ii) that $Pf$ satisfies the homogeneous Weinstein equation while $Q f$ satisfies the inhomogeneous Weinstein equation with parameter $-k$.\\
Applying now the Laplacian to $h:=\frac{f e_n}{x_n^k}$, yields
 \begin{eqnarray*}
  \bigtriangleup \left( \frac{f e_n}{x_n^k} \right)  &=& k(k+1) \frac{f e_n}{x_n^{k+2}} + 2 (-k) \frac{1}{x_n^{k+1}}(\frac{\partial}{\partial x_n} f) e_n + \frac{\bigtriangleup f e_n}{x_n^k}\\
&\stackrel{\scriptscriptstyle Weinstein}{=}& k(k+1) \frac{f e_n}{x_n^{k+2}} + 2 (-k) \frac{1}{x_n^{k+1}}(\frac{\partial}{\partial x_n} f) e_n \\
&&+ \frac{k}{x_n^{k+1}}(\frac{\partial}{\partial x_n} f) e_n  - k \frac{(Qf)e_n}{x_n^{k+2}}e_n\\
&=& k(k+1) \frac{f e_n}{x_n^{k+2}} -k \frac{1}{x_n^{k+1}}(\frac{\partial}{\partial x_n} f) e_n - k \frac{(Qf)e_n}{x_n^{k+2}} e_n \\
&=& k \frac{f e_n}{x_n^{k+2}} -k \frac{1}{x_n}\frac{\partial}{\partial x_n} (\frac{f e_n}{x_n^k}) - k \frac{(Qf)e_n}{x_n^{k+2}} e_n\\
&=& k\frac{(Pf)e_n}{x_n^{k+2}} - k \frac{1}{x_n}\frac{\partial}{\partial x_n} (\frac{f e_n}{x_n^k})\\
&=& k Q(\frac{f e_n}{x_n^k}) e_n - k \frac{1}{x_n}\frac{\partial}{\partial x_n} (\frac{f e_n}{x_n^k})
 \end{eqnarray*}
So $h$ satisfies the Weinstein equation for $-k$ and thus $ \bigtriangleup^{\frac{-k+2}{2}} h = 0$ with Theorem \ref{satz_eins}. As the kernel of any power of $\bigtriangleup$ is invariant under multiplication with $e_n$ from the right, we also obtain $ \bigtriangleup^{\frac{-k+2}{2}} g = 0$.

\end{proof}
Now we can readily derive that 
\begin{corollary}
 Let $k$ be an even integer. Assume that $f$ is sufficiently many times real differentiable function with
\[
  D \bigtriangleup^{\frac{k}{2}} f = 0.
\]
Then the function $g:=\frac{f}{x_n^k}$ satisfies 
\[
  D \bigtriangleup^{\frac{-k+2}{2}} g = 0.
\]
\end{corollary}
\begin{proof}
Suppose that $k \in 2 \mathbb{Z}$. Let $f \in Ker D \Delta^{k/2}$.
Then we automatically have that $f \in $ Ker $\Delta^{(k+2)/2}$. Using Theorem \ref{satz_zwei}, we know that $g:=f/x_n^k$ is an element from  Ker $\Delta^{(-k+2)/2}$. But then of course we also have 
$g \in $ Ker $D \Delta^{(-k+2)/2}$. So $g$ is $(-k+2)$-holomorphic Cliffordian.
\end{proof}

Therefore, it is sufficient to restrict to either non-negative or non-positive values of $k$ in all that follows. Notice that we only look at even values for $k$. Finally, the set of $k$-holomorphic Cliffordian functions is invariant under the action of $SAV(\mathbb{R} \oplus \mathbb{R}^{n})$. 
In view of the invariance of the Weinstein equation under M\"obius transformations, cf. \cite{AL}, one can directly adapt the proof of the $k$-hypermonogenic case given in \cite{slnewnew} to establish:
\begin{theorem} Let $k \in 2{\mathbb{Z}}$. Suppose that $M \in SAV(\mathbb{R} \oplus \mathbb{R}^n)$. Let    
 $y=M<x>=(ax+b)(cx+d)^{-1}$ be the image of a point $x$ under such a M\"obius transformation.  Then such a function $f(y)$ that is $k$-holomorphic Cliffordian in the variable $y$ is transformed to a function
\begin{equation}
\label{kinv1}
F(x):=\frac{\overline{cx+d}}{\|cx+d\|^{n+1-k}}f(M<x>)
\end{equation}
which turns out to be $k$-holomorphic Cliffordian in the variable $x$.
\end{theorem}

Due to the equivalence described in Theorem 1(ii) between the solutions of the Weinstein equation and the null-solutions to $\Delta^{(k+2)/2}$, the invariance property of the Weinstein equation under M\"obius transformations is automatically inherited for the whole class of functions satisfying $\Delta^{(k+2)/2} f = 0$.  For the subclass of $k$-holomorphic Cliffordian functions, the canonical automorphy factor is the $k$-holomorphic Cliffordian kernel 
$G_k(x)=\frac{\overline{x}}{\|x\|^{n+1-k}}$. This is the same expression as the kernel for the subset of $k$-hypermonogenic functions which is contained in the set of $k$-holomorphic Cliffordian functions. 

\par\medskip\par

This invariance property allows us to introduce $k$-holomorphic Cliffordian automorphic forms.

\section{$k$-holomorphic Cliffordian automorphic forms}

\subsection{Basic properties and examples}
Now we have all tools in hand to introduce 
\begin{definition}
Let $p \le n-1$ and suppose that $k \in 2 {\mathbb{Z}}$. 
A left $k$-holomorphic Cliffordian function $f: H^{+}(\mathbb{R} \oplus\mathbb{R}^n) \rightarrow Cl_n$ is called a left $k$-holomorphic Cliffordian  automorphic form on $\Gamma_{p}({\cal{I}})[N]$ if for all $x \in H^{+}(\mathbb{R} \oplus\mathbb{R}^n)$  
\begin{equation}
\label{stern}
f(x)=\frac{\overline{cx+d}}{\|cx+d\|^{n+1-k}} f(M<x>) 
\end{equation}
for all $M \in \Gamma_{p}({\cal{I}})[N]$. 
\end{definition}
In the case $k=0$ we re-obtain the class of left monogenic automorphic forms discussed in \cite{KraHabil}. All $k$-hypermonogenic automorphic forms discussed in \cite{CKR2007} are included in the set of $k$-holomorphic Cliffordian automorphic forms. 
As the following proposition shows, there is a direct relation between $k$-holomorphic Cliffordian automorphic forms and  $-k+2$-holomorphic Cliffordian automorphic forms.  
\begin{proposition}\label{-k}
Suppose that $k\in 2\mathbb{Z}$ and that $p$ is a positive integer with $p < n$. If $f:H^{+}(\mathbb{R} \oplus\mathbb{R}^n) \rightarrow Cl_n$ is a $k$-holomorphic Cliffordian automorphic form on $\Gamma_{p}({\cal{I}})[N]$ satisfying (\ref{stern}), then the function $g:H^{+}(\mathbb{R} \oplus\mathbb{R}^n) \rightarrow Cl_n$ defined by $g(x) := \frac{f(x)}{x_n^k}$ is $-k+2$ holomorphic Cliffordian and satisfies 
\begin{equation}
g(x) =  \frac{\overline{cx+d}}{\|cx+d\|^{n+1+k}} g(M<x>)
\end{equation}
for all $M \in \Gamma_{p}({\cal{I}})[N]$. 
\end{proposition}
\begin{proof}
Suppose that $f$ is a $k$-holomorphic Cliffordian automorphic form satisfying (\ref{stern}). Then  
$g(x) := \frac{f(x)}{x_n^k}$ 
is $-k+2$-holomorphic Cliffordian and satisfies the transformation law 
\begin{eqnarray*}
g(M<x>) &=& \frac{1}{(M<x>)_n^k} f(M<x>)\\
&=& \frac{\|cx+d\|^{2k}}{x_n^k} (\|cx+d\|^{n+1-k}) (\overline{cx+d})^{-1} f(x)\\
&=& \|cx+d\|^{n+1+k}(\overline{cx+d})^{-1} g(x).
\end{eqnarray*}
\end{proof}
Next we want to construct non-trivial examples of $k$-holomorphic Cliffordian automorphic forms. Before we proceed in that way we need to recall an important algebraic fact.
For each ${\cal{I}}$ there exists a minimal positive integer $N_0({\cal{I}})$ such that neither the negative identity matrix $-I$ nor the other diagonal matrices of the form $\left(\begin{array}{cc} e_A^* & 0 \\ 0 & e_A^{-1} \end{array} \right)$ where $A \subseteq P(1,\ldots,p)$ are not any longer included in all principal congruence subgroups $\Gamma_{p}({\cal{I}})[N]$ with $N \ge N_0({\cal{I}})$. In the case where ${\cal{I}}$ is the standard $\mathbb{Z}$-order ${\cal{O}}_p$, we have $N_0({\cal{O}}_p) = 3$, see \cite{KraHabil}. For all $N < N_0({\cal{I}})$, only the zero function satisfies (\ref{stern}). However, for all $N \ge N_0({\cal{I}})$, one can construct non-trivial $k$-hypermonogenic automorphic forms that have the transformation behavior (\ref{stern}).
\par\medskip\par

For non-positive even integers $k$, the simplest non-trivial examples of $k$-holomorphic Cliffordian automorphic forms on the groups $\Gamma_p({\cal{I}})[N]$ with $N \ge N_0({\cal{I}})$ are the $k$-hypermonogenic generalized Eisenstein series, introduced in \cite{CKR2007}.

\begin{definition}
Let $N \ge N_0({\cal{I}})$. For $p < n$ and even integer $k$ with $k < n-p-1$ the $k$-hypermonogenic Eisenstein series on the group $\Gamma_{n-1}({\cal{I}})[N]$ acting on $H^{+}(\mathbb{R} \oplus\mathbb{R}^n)$ are defined by
\begin{equation}\label{khypereisen}
\varepsilon_{k,p,N}(x) := \sum\limits_{M:{\cal{T}}_{p}({\cal{I}})[N] \backslash \Gamma_{p}({\cal{I}})[N]} \frac{\overline{cx+d}}{\|cx+d\|^{n+1-k}}.
\end{equation}
\end{definition}
These series converge for $k < n-p-1$ absolutely and uniformly on each compact subset of $H^{+}(\mathbb{R} \oplus\mathbb{R}^n)$. A majorant is 
$$
\sum\limits_{M: {\cal{T}}_p({\cal{I}})[N] \backslash \Gamma_p({\cal{I}})[N]} \frac{1}{\|c e_n +d \|^{\alpha}}
$$
whose absolute convergent abscissa is $\alpha > p+1$, cf. for example \cite{KraHabil,EGM90}. For $p=n-1$, this majorant converges absolutely for all $k <-1$.  

\par\medskip\par

The non-vanishing behavior for $N \ge N_0({\cal{I}})$ can easily be established by considering the limit $\lim_{x_n \rightarrow +\infty} \varepsilon_{k,p,M}(x_n e_n)$ which equals $+1$ in these cases,~cf.~\cite{KraBMS}.  

\par\medskip\par 

{\bf Remarks.}  

In the case $k=0$ and ${\cal{I}} = {\cal{O}}_p$ the series (\ref{khypereisen}) coincide with the monogenic Eisenstein series considered in \cite{KraHabil,KraBMS} in the cases $p < n-2$.  

\par\medskip\par

By adapting the Hecke trick from \cite{Freitag} one can also introduce Eisenstein series of lower weight. This is shown in Section 4 of \cite{paper4}. In particular, for $N \ge N_0({\cal{I}})$ the series 
\begin{equation}\label{0hypereisen}
\varepsilon_{0,N-1,M}(x) := \lim\limits_{s \rightarrow 0^+}\sum\limits_{M:{\cal{T}}_{n-1}({\cal{I}})[N] \backslash \Gamma_{n-1}({\cal{I}})[N]} \Bigg(\frac{x_n}{\|cx+d\|^2}\Bigg)^s \frac{\overline{cx+d}}{\|cx+d\|^{n+1}}.
\end{equation}
defines a well-defined non-vanishing left monogenic Eisenstein series on the groups $\Gamma_{n-1}({\cal{I}})[N]$ in  upper half-space variable $x$, cf.~\cite{paper4}~Section~4. In \cite{paper4}~Section~4 this is done for the particular case ${\cal{I}} = {\cal{O}}_p$. However, the transition to the context of more general orders ${\cal{I}}$ follows identically along the same lines.   

\par\medskip\par

In view of Proposition~\ref{-k}, we can directly construct non-vanishing $j+2$-holomorphic Cliffordian Eisenstein series for positive $j$ from the $k$-hypermonogenic Eisenstein series of negative $k$, simply by forming
$$
E_{-k,p,N} (x) := \frac{\varepsilon_{k,p,N}(x) e_n}{x_n^k} = \sum\limits_{M:{\cal{T}}_{p}({\cal{I}})[N] \backslash \Gamma_{p}({\cal{I}})[N]} \frac{(\overline{cx+d})e_n}{\|cx+d\|^{n+1-k}x_n^k}.
$$
These series then satisfy the transformation law
$$
E_{-k,p,N} (x) := \frac{\overline{cx+d}}{\|cx+d\|^{n+1+k}} E_{-k,p,N}(M<x>)
$$    
for all $M \in \Gamma_p({\cal{I}})[N]$. The functions $E_{-k,p,N} (x)$ are in Ker $D \Delta^{(j+2)/2}$ where $j=-k > 0$. Since they are even $j$-hypermonogenic, we even have that they are $j$-holomorphic Cliffordian. This is a stronger statement than being $j+2$-holomorphic Cliffordian. 

It is clear that the $j$-holomorphic Cliffordian Eisenstein series 
$$
E_{j,p,N}(x) := \sum\limits_{M:{\cal{T}}_{p}({\cal{I}})[N] \backslash \Gamma_{p}({\cal{I}})[N]} \frac{x_n^j(\overline{cx+d})e_n}{\|cx+d\|^{n+1+j}}
$$
converge for every $j > 1$. They provide us with the simplest examples of $j$-holomorphic Cliffordian automorphic forms of weight $(n-j)$ on the full group $\Gamma_{n-1}({\cal{I}})[N]$ with $N \ge 3$ for positive $j>1$. In the $k$-holomorphic Cliffordian setting, we may omit the factor $e_n$ on the right-hand side. 

\par\medskip\par
{\bf Remark}: For the case $k=n-1$ with $n$ odd G. Laville and I. Ramadanoff considered $k$-holomorphic Cliffordian generalizations of the cotangent function, the Weierstra{\ss} $\zeta$-function and the Weierstra{\ss} $\wp$-function. These are the simplest examples of $k$-holomorphic Cliffordian automorphic forms for discrete translation groups. See \cite{LR2,LR3}. Analogous constructions for polymonogenic functions (i.e. functions in the kernel of arbitrary iterates of the Dirac operator) are presented in \cite{KraHabil}.

\subsection{Relation to Maa{\ss} forms}

As one directly sees, the theory of monogenic automorphic forms fits as a special case within the general framework of  $k$-holomorphic Cliffordian automorphic forms. We can say more:

\par\medskip\par

Suppose that $f = Pf + Qf e_n$ is a $k$-holomorphic Cliffordian function, where $k$ is some arbitrary fixed integer number. Then, as stated in Theorem~1 (ii) the real components of $Pf$ are $k$-hyperbolic harmonic functions, i.e. solutions to 
\begin{equation}
\label{khypharm}
x_n \Delta u- k \frac{\partial u}{\partial x_n} = 0.
\end{equation}
Remember that this is also true for the particular subset of $k$-hypermonogenic functions. 
Following \cite{slnew} and others, the solutions to~(\ref{khypharm}) are also quasi-invariant under M\"obius transformations that act on upper half-space: If $f$ is a solution to (\ref{khypharm}), then 
\begin{equation}
\label{harminv}
F(x) = \frac{1}{\|cx+d\|^{n-k-1}} f(M<x>)
\end{equation}
is $k$-hyperbolic harmonic, too. 
Notice that in the particular case $k=n-1$ the correction factor disappears. This property attributes a special role to the function  class of $(n-1)$-holomorphic Cliffordian functions, whence we are in the function class considered in \cite{LR} (supposed we are in an even dimensional space).  

\par\medskip\par

As mentioned earlier in \cite{CKR2007}, the solutions of (\ref{khypharm}) in turn are directly related to the Maa{\ss} wave equation. Following for example \cite{Leutw87}, if $u$ is a solution of (\ref{khypharm}) then  $g(x) = x_n^{-(1-n+k)/2}u(x)$ is a solution of
\begin{equation}
\label{maass}
\Delta g - \frac{n-1}{x_n} \frac{\partial g}{\partial x_n} + \lambda \frac{g}{x_n^2} = 0
\end{equation}
where $\lambda = \frac{1}{4}(n^2 -(k+1)^2)$. 
The solutions to (\ref{maass}) have the property that they are directly preserved by all M\"obius transformations that  act on upper half-space. Each solution $g$ is an eigensolution to the Laplace-Beltrami operator (\ref{LB}) associated to the  fixed eigenvalue $-\frac{1}{4}(n^2-(k+1)^2).$ 

\par\medskip\par

Let $p < n$. Now suppose that $f: H^{+}(\mathbb{R} \oplus\mathbb{R}^n) \rightarrow Cl_n$ is a $k$-holomorphic Cliffordian automorphic form on $\Gamma_{p}({\cal{I}})[N]$ of weight $(n-k)$, satisfying the transformation law $f(x) = \frac{\overline{cx+d}}{\|cx+d\|^{n+1-k}} f(M<x>)$ for all $M \in \Gamma_{p}({\cal{I}})[N]$.

Since $(M<x>)_n = \frac{x_n}{\|cx+d\|^2}$, the function $g(x) = x_n^{-(1-n+k)/2}f(x)$ thus satisfies for all $M \in \Gamma_{p}({\cal{I}})[N]$:
\begin{eqnarray*}
g(M<x>) &=& {(M<x>)_n}^{-\frac{1-n+k}{2}} f(M<x>) \\
&=& \Bigg(\frac{x_n}{\|cx+d\|^2}\Bigg)^{-\frac{1-n+k}{2}} \|cx+d\|^{n+1-k} (\overline{cx+d})^{-1} f(x)\\
&=& x_n^{-\frac{1-n+k}{2}} \frac{1}{\|cx+d\|^{n-k+1}} \|cx+d\|^{n+1-k} (\overline{cx+d})^{-1} x_n^{\frac{1-n+k}{2}} g(x).
\end{eqnarray*}
Hence, 
$$
g(x) = \frac{\overline{cx+d}}{\|cx+d\|^2} g(M<x>). 
$$
Unfortunately, if $f = Pf + Qf e_n$ is a $k$-holomorphic Cliffordian automorphic form with respect to $\Gamma_{p}({\cal{I}})[N]$, then $Pf$ is in general not an automorphic form with respect to the full group $\Gamma_p[N]$. However, only the components of the $P$-part of $f$ satisfy the homogeneous Weinstein equation (\ref{khypharm}). The associated function $g(x)= x_n^{-(1-n+k)/2}f(x)$ is exactly a $\Gamma_p({\cal{I}})[N]$-invariant eigenfunction to the Laplace-Beltrami operator for the eigenvalue  $-\frac{1}{4}(n^2-(k+1)^2)$, if $Q f= 0$ or if $k=0$. If $k=0$, then also the $Q$-part of $f$ satisfies the homogeneous Weinstein equation, due to the cancellation of the term on the right-hand side. The latter case corresponds exactly to the monogenic case.  We cite from \cite{CKR2007}:
\begin{theorem}\label{relation}
Suppose that $f$ is a left monogenic automorphic form on $\Gamma_{p}({\cal{I}})[N]$ of weight $n$, satisfying 
$f(x) = \frac{\overline{cx+d}}{\|cx+d\|^{n+1}} f(M<x>)$ for all $M \in \Gamma_{p}({\cal{I}})[N]$.  Then $g(x) = x_n^{-(1-n)/2}f(x)$ is a quasi-$\Gamma_p({\cal{I}})[N]$-invariant Maa{\ss} wave form associated to the fixed eigenvalue $-\frac{1}{4}(n^2-(k+1)^2)$ and has the $-1$-weight automorphy factor $\frac{\overline{cx+d}}{\|cx+d\|^2}$.   
\end{theorem}

{\bf Remarks:}  
If $F$ is a $k$-hyperbolic harmonic automorphic form that satisfies $F(x) = \|cx+d\|^{-(n-k-1)} F(M<x>)$, then $G(x)= x_n^{-(1-n+k)/2}f(x)$ turns out to be an eigenfunction to the Laplace-Beltrami operator associated to the eigenvalue $-\frac{1}{4}(n^2-(k+1)^2)$. The function $G$ then is totally invariant under the group action of $\Gamma_{p}({\cal{I}})[N]$, i.e. $G(x) = G(M<x>)$ for all $M \in \Gamma_{p}({\cal{I}})[N]$.

Theorem~\ref{relation} exhibits an explicit connection between the class of  $k$-holomorphic Cliffordian automorphic forms (in which the monogenic automorphic forms play a particular role) and the particular family of non-analytic automorphic forms on the Ahlfors-Vahlen group considered for instance by A. Krieg, J. Elstrodt et al. and others (see for example \cite{Kri88,EGM90,Kri03}). The non-analytic automorphic forms considered in \cite{Kri88,EGM90,Kri03}) are scalar-valued eigenfunctions of the Laplace-Beltrami operator 
associated to a special continuous spectrum of eigenvalues described in \cite{EGM90}. They are all totally invariant under the group action. 

\par\medskip\par

Notice that the Eisenstein- and Poincar\'e series considered in this paper are Clifford algebra-valued in general. They are associated to a fixed eigenvalue. In the case where they are not scalar valued and where $k \neq 0$, their $Q$-part satisfies the inhomogeneous Weinstein equation with parameter $-k$ instead of the hyperbolic harmonic equation. In these cases only their $P$-parts are $k$-hyperbolic harmonic.   

\section{Cusp forms and orthogonality relations}

All the classical complex-analytic Eisenstein- and Poincar\'e series on congruence subgroups of finite index in the modular group  $SL(2,\mathbb{Z})$ have the property that they are orthogonal to each other with respect to the Petersson  inner  product. The attached space of the complex-analytic Eisenstein series is exactly the orthogonal complement of the space of the complex-analytic cusp forms on these groups. Cusp forms are those automorphic forms that vanish at each cusp of the group. In turn, these are exactly those automorphic forms which have a vanishing first Fourier coefficient. 
 
The complex-analytic Poincar\'e series are special examples of cusp forms. For details we refer for instance to  \cite{Freitag}. An analogous result has been established for the scalar valued non-analytic Maa{\ss} wave forms on arithmetic subgroups of the Ahlfors-Vahlen group in \cite{EGM90}.

\par\medskip\par

In this section we want to establish a similar result within the context of $k$-holomorphic Cliffordian functions. Notice that an extension of the classical series ansatz
$$
P(x,w) = \sum\limits_{M \in SL(2,\mathbb{Z})} (cx+d)^{-k}(w+M<x>)^{-k}
$$
does not lead to a sensible construction of $k$-hypermonogenic cusp forms, because the translation with $w$ (having a non-zero $e_n$-component) destroys $k$-hypermonogenicity. This was an important reason for us to consider in this paper the larger class of $k$-holomorphic Cliffordian functions. 

\par\medskip\par

In all that follows we restrict to consider those particular orders ${\cal{I}}$ which have the property that $+\infty$ is the only cusp of ${\cal{I}}$ and that the quotient $H^+(\mathbb{R} \oplus{\mathbb{R}^n}) \cup \{ {\rm cusps\;of\;}\Gamma_{n-1}({\cal{I}})[N] \}/\Gamma_{n-1}({\cal{I}})[N]$ is compact. In the cases $n < 4$ the orders ${\cal{I}} = {\cal{O}}_{n-1}$ provide us with some simple examples having this property.    

\subsection{Cusp forms}
To start we introduce:
\begin{definition}\label{moncusp} ($k$-holomorphic Cliffordian cusp form)

For even integers $k \le 0$ a left $k$-holomorphic Cliffordian cusp form on $\Gamma_{n-1}({\cal{I}})[N]$ is a left $k$-holomorphic Cliffordian automorphic form on $\Gamma_{n-1}({\cal{I}})[N]$ that satisfies additionally
\begin{equation}\label{cusp}
\lim\limits_{x_n \rightarrow +\infty} x_n^{-k} \frac{\overline{cx_n e_n + d}}{\|cx_n e_n + d\|^{n+1-k}} f(M<x_n e_n>) = 0 \quad \quad \forall M \in \Gamma_{n-1}({\cal{I}})[N].
\end{equation}
For positive even integers $k$, the factor $x_n^{-k}$ is omitted.
\end{definition}
{\bf Remark}: In view of the quasi-invariance of $f$ under $\Gamma_{n-1}({\cal{I}})[N]$, we deal with a cusp form if (\ref{cusp}) is satisfied for all matrices $R \in {\cal{R}}$, where ${\cal{R}}$ is any system of representatives of right cosets in $\Gamma_{n-1}({\cal{I}})$ modulo $\Gamma_{n-1}({\cal{I}})[N]$. More precisely: This relation must hold for any arbitrary system of representatives. 
Since $\Gamma_{n-1}({\cal{I}})[N]$ is a subgroup of finite index in $\Gamma_{n-1}({\cal{I}})$, any such system  ${\cal{R}}$ has a finite cardinality.  

\par\medskip\par
For simplicity let us restrict to the case $k < 0$. 
Concrete examples of left $k$-holomorphic Cliffordian cusp forms on $\Gamma_{n-1}({\cal{I}})[N]$, $N \ge N_0({\cal{I}})$ arise from the following family of Poincar\'e series defined by   
\begin{equation}\label{khp}
P_{k,n-1,N}(x,w) := \sum\limits_{M \in \Gamma_{n-1}({\cal{I}})[N]} \frac{\overline{cx+d}}{\|cx+d\|^{n+1-k}} \frac{\overline{w+M<x>}}{\|w+M<x>\|^{n+1-k}},
\end{equation}
whenever the quotient $H^+(\mathbb{R} \oplus {\mathbb{R}^n}) \cup \{ {\rm cusps\;of\;} \Gamma_{n-1}({\cal{I}})[N] \}/\Gamma_{n-1}{\cal{I}}[N]$ is compact.

This series is majorized by the series 
\begin{equation}\label{khpm}
\sum\limits_{M \in \Gamma_{n-1}({\cal{I}})[N]} \frac{C}{\|cx+d\|^{n-k}\|w+M<x>\|^{n-k}},
\end{equation} 
for some properly chosen real $C > 0$. 
The absolute convergence abscissa of (\ref{khpm}) on $H^{+}(\mathbb{R} \oplus \mathbb{R}^n) \times H^{+}(\mathbb{R} \oplus \mathbb{R}^n)$ is $k <-1$.  This can be shown be applying the usual compactification arguments, see Lemma~1 from \cite{ck} for details.  

\par\medskip\par

In view of Corollary 1, for even positive $k > 1$ we can introduce $k+2$-holomorphic Cliffordian Poincar\'e series by 
\begin{equation}\label{khp+}
P_{k,n-1,N}(x,w) := x_n^k \sum\limits_{M \in \Gamma_{n-1}({\cal{I}})[N]} \frac{\overline{cx+d}}{\|cx+d\|^{n+1+k}} \frac{\overline{w+M<x>}}{\|w+M<x>\|^{n+1+k}}. 
\end{equation}

\par\medskip\par

By a direct calculation one can show that these series really vanish at each cusp of the group. Notice that this series actually is not $k$-hypermonogenic. This is due to the application of the shift of $w$ in the expression $w+M<x>$. Since $w$ as a non-vanishing $e_n$-component, $k$-hypermonogenicity is not preserved. However, the function remains $k$-holomorphic Cliffordian.
\par\medskip\par 
An important task is to establish  
\begin{theorem} (Non-vanishing argument)\\
Let $k < -1$ be an even integer and $N \ge N_0({\cal{I}})$. Then the series (\ref{khp}) do not vanish identically for all $k$.  
\end{theorem}
\begin{proof} 
First, we choose a fixed $w$ and $x$ from $H^{+}(\mathbb{R} \oplus\mathbb{R}^n)$. By applying the usual compactification argument as one uses in the convergence proof, we can conclude that  
$$
\mu:=\min_{M \in \Gamma_{n-1}({\cal{I}})[N]}\{\|cx+d\| \|w+M<x>\|\} > 0.
$$
The miminum can be attained by several matrices, say by $M_1,\ldots,M_k$. Since $M$ is a discrete group, the number of these matrices is finite. Otherwise we would have an accumulation point. All the other matrices $M'$ satisfy 
$$
T_{M'}(x):= \|cx+d\| \|w+M'<x>\| > \mu.
$$
We thus have 
\begin{equation}\label{Matrixminimum}
|{T_M}_i| < |T_{M'}| - \varepsilon
\end{equation}
for all $M_1,\ldots,M_k$ and $M' \in \Gamma_{n-1}({\cal{I}})[N]$.
In view of the continuity of the expressions $T_M(x)$, the estimate remains valid in a sufficiently small neighborhood of $w$, say in $W$. Let us now split the Poincar\'e series (\ref{khp}) as follows:
\begin{eqnarray*}
P_{k,n-1,N}(x,w) &=& \sum\limits_{i=1,\ldots,k}   \frac{\overline{cx+d}}{\|cx+d\|^{n+1-k}} \frac{\overline{w+M_i<x>}}{\|w+M_i<x>\|^{n+1-k}} 
\end{eqnarray*}
\begin{eqnarray*}
& &+ 
\sum\limits_{M \in \Gamma_{n-1}({\cal{I}})[N] \backslash \{M_1,\ldots,M_k\}} \frac{\overline{cx+d}}{\|cx+d\|^{n+1-k}} \frac{\overline{w+M<x>}}{\|w+M<x>\|^{n+1-k}}.
\end{eqnarray*}
Let us denote the first part of this sum by $P_1$ and the second part by $P_2$. The first sum is a finite sum. It hence  possesses an analytic extension in $w$ up to point singularities to the lower half-space. 

Under the condition $N \ge N_0({\cal{I}})$, the finite series $P_1$ does not vanish identically, since the expression $\frac{1}{\|w+M<x>\|}$ has an isolated singularity at $w=-M<x>$. The condition $N \ge N_0({\cal{I}})$ is crucially applied for this conclusion. 
For all $N \ge N_0({\cal{I}})$ the expression $M<x>$ is reached by $M$ only.   

Next we observe that $P_1$ grows asymptotically like $\Omega(\frac{1}{\mu^{n-k}})$ in $W$. 
As to $P_2$, we can estimate by using (\ref{Matrixminimum}) that
\begin{eqnarray*}
 &&  \max\limits_{M \in \Gamma_{n-1}({\cal{I}})[N] \backslash \{M_1,\ldots,M_k\}} \left|\frac{\overline{cx+d}}{\|cx+d\|^{n+1-k}} \frac{\overline{w+M<x>}}{\|w+M<x>\|^{n+1-k}} \right| \\
&=& \max\limits_{M \in \Gamma_{n-1}({\cal{I}})[N] \backslash \{M_1,\ldots,M_k\}} \frac{1}{\|T_M\|^{n-k}} \\
&<&  \frac{1}{|\mu + \varepsilon|^{n-k}},
\end{eqnarray*}
This tells us that $P_2$ grows like $\mathcal{O} (\frac{1}{ ( \mu + \varepsilon )^{n-k} } )$.
\par\medskip\par

Suppose now that the complete Poincar\'e series $P_{k,n-1,N}(x,w) $ vanishes identically for all $w$ and for all $k<-1$.  Since $P_1$ and $P_2$ have a different asymptotic growth behavior in $W$, it follows that both the partial series $P_1$ and $P_2$ vanish identically. This is a contradiction. 

\end{proof}

{\bf Remark.} Applying Corollary 1 to the preceding statement allows us to establish the existence of $k$-holomorphic Cliffordian cusp forms for infinitely many positive even integers $k>3$. 

\subsection{The Fourier expansion of $k$-holomorphic Cliffordian automorphic forms} 

In this section we show that $k$-holomorphic Cliffordian cusp forms are exactly those $k$-holomorphic Cliffordian automorphic forms whose Fourier image of zero-frequency vanishes identically. This is a nice analogue to the classical cases. 

\par\medskip\par

To proceed in this direction we first determine the Fourier series development of a general $k$-holomorphic Cliffordian function on upper half-space. To do so we recall that the $P$-part of a $k$-holomorphic Cliffordian function satisfies the homogeneous Weinstein equation and the $Q$-part the inhomogeneous one associated with the parameter $-k$. This is a consequence of the fact that Ker $D\Delta^{k/2} \subset$ Ker $\Delta^{(k+2)/2}$.  
We write $f$ in the form $f=Pf+Qf e_n$. We obtain 
\begin{eqnarray*}
\Bigg(\Delta-\frac{k}{x_n} \frac{\partial }{\partial x_n}\Bigg) f & = & \Bigg(\Delta-\frac{k}{x_n} \frac{\partial }{\partial x_n}\Bigg) (Pf + Qf e_n) \\
&=& \Bigg(\Delta-\frac{k}{x_n} \frac{\partial }{\partial x_n}\Bigg) Pf + \Bigg(\Delta-\frac{k}{x_n} \frac{\partial }{\partial x_n}\Bigg) (Qf) e_n \\
&=& 0 - \frac{k}{x_n} Qf e_n.
\end{eqnarray*}
Summarizing, we obtained that every $k$-holomorphic Cliffordian function satisfies the equation 
$$
\Bigg(\Delta-\frac{k}{x_n} \frac{\partial }{\partial x_n}\Bigg) f = - \frac{k}{x_n} Qf e_n.
$$
When applying the Fourier transform in the first $n$ variables $x_0,x_1,\cdots,x_{n-1}$, we obtain (using the decomposition of $f$ into its $P$- and $Q$-part): 
\[
\sum_{r=0}^{n-1} (i \omega_r)^2 \alpha + \frac{\partial^2 }{\partial x_n^2} \alpha - \frac{k}{x_n}
\frac{\partial }{\partial x_n} \alpha - e_n k \frac{(Q\alpha)'}{x_n^2} = 0.
\]
Here $\underline{\omega} = \omega_0 + \omega_1 e_1 + \cdots + \omega_{n-1} e_{n-1}$ is the frequency vector in the domain of the Fourier image. 
The latter equation corresponds to 
\[
-\|\underline{\omega}\|^2 \alpha + \frac{\partial^2 }{\partial x_n^2} \alpha - \frac{k}{x_n} \frac{\partial }{\partial x_n} \alpha - e_n k
\frac{(Q\alpha)'}{x_n^2} = 0.
\]
Now we put 
\[
 \alpha = P \alpha - e_n (Q\alpha)'
\]
and obtain the system 
\begin{eqnarray*}
- \|\underline{\omega}\|^2 \alpha &=& -\|\underline{\omega}\|^2 P \alpha + e_n \|\underline{\omega}\|^2 (Q\alpha)' \\
 \frac{\partial^2 }{\partial x_n^2} \alpha &=& \frac{\partial^2}{\partial x_n^2} P \alpha - e_n \frac{\partial^2 }{\partial x_n^2} (Q\alpha)' \\
- \frac{k}{x_n} \frac{\partial }{\partial x_n} \alpha &=& - \frac{k}{x_n} \frac{\partial }{\partial x_n} P \alpha + e_n
\frac{k}{x_n}\frac{\partial }{\partial x_n}(Q\alpha)' \\
-e_n k \frac{(Q\alpha)'}{x_n^2} &=& -e_n  \frac{k}{x_n^2} \left( Q\alpha \right)'. 
\end{eqnarray*}
This leads to 
\begin{equation}\label{geins}
- \|\underline{\omega}\|^2 P \alpha + \frac{\partial^2 }{\partial x_n^2} P \alpha - \frac{k}{x_n}\frac{\partial }{\partial x_n} P \alpha  = 0
\end{equation}
\begin{equation}\label{gzwei}
 \|\underline{\omega}\|^2 (Q\alpha)' - \frac{\partial^2 }{\partial x_n^2} (Q\alpha)' + 
\frac{k}{x_n}\frac{\partial }{\partial x_n}(Q\alpha)' -\frac{k}{x_n^2} \left( Q\alpha \right)' = 0
\end{equation}
This system is equivalent to 
\begin{eqnarray*}
 \frac{\partial^2 }{\partial x_n^2} P \alpha &=& \|\underline{\omega}\|^2 P \alpha +  \frac{k}{x_n} \frac{\partial }{\partial x_n} P \alpha  \\
  \frac{\partial^2 }{\partial x_n^2} (Q\alpha)' &=& \|\underline{\omega}\|^2 (Q\alpha)' + 
\frac{k}{x_n} \frac{\partial }{\partial x_n}(Q\alpha)' - \frac{k}{x_n^2} \left( Q\alpha \right)'.
\end{eqnarray*}
The requirement that the limit as $x_n\to\infty$ must vanish, determines that $P\alpha$ must be proportional to $x_n^{(k+1)/2}K_{(k+1)/2}(\|\underline{\omega}\|x_n)$ and $(Q\alpha)'$ to the expression  $x_n^{(k+1)/2}K_{(k-1)/2}(\|\underline{\omega}\|x_n)$. 
These two proportionality coefficients are furthermore independent of each other. The zero frequency solutions are clearly $1$ and $x_n^k$ for $P$ and $x_n$ and $x_n^{k+1}$ for $Q$, all with independent coefficients. We can now establish 

\begin{theorem}\label{kfnew}
Let $f: H^{+}(\mathbb{R} \oplus \mathbb{R}^n) \rightarrow Cl_n$ be a $k$-holomorphic Cliffordian automorphic form on $\Gamma_{n-1}({\cal{I}})[N]$. Further, let ${\cal{R}}$ be a system of representatives in $\Gamma_{n-1}({\cal{I}})$ modulo $\Gamma_{n-1}({\cal{I}})[N]$.\\ Then for every $R \in {\cal{R}}$ the function $F(x;R):= \frac{\overline{(c_R x+d_R)}}{\|c_R x+d_R\|^{n+1-k}}f(R<x>)$ has a particular Fourier series representation of the form 
{\small{
\begin{eqnarray}
\label{khfourier}
F(x;R) &=&  a(\underline{0};R) + \alpha(\underline{0};R) x_n^k + b(\underline{0};R) e_n x_n + \beta(\underline{0};R) e_n x_n^{k+1} \\
&+&\!\!\!\!\!  \sum_{\underline{\omega} \in \Lambda_{n-1}^{*}[N] \backslash \{\underline{0}\}} \!\!
x_n^{\frac{k+1}{2}} \Bigg[K_{\frac{k+1}{2}}(2 \pi \|\underline{\omega}\| x_n) \alpha(\underline{\omega};R) + K_{\frac{k-1}{2}}(2 \pi \|\underline{\omega}\| x_n) \beta(\underline{\omega};R) e_n \Bigg] \nonumber \\
& & \quad\quad\quad\quad\quad \times 
 e^{2 \pi i <\underline{\omega},\underline{x}>}.\nonumber
\end{eqnarray}}}
Here, $a(0,R),b(0;R)$ and all $\alpha(\underline{\omega};R)$ and $\beta(\underline{\omega};R)$ are well-defined Clifford numbers from the Clifford sub algebra $Cl_{n-1}$. 
\end{theorem}
Let us now compare this result with the special $k$-hypermonogenic setting. From \cite{CKR2007} we recall that in the $k$-hypermonogenic case, one has the following particular form of Fourier expansion:  
\begin{theorem}\label{kf}
Let $f: H^{+}(\mathbb{R} \oplus \mathbb{R}^n) \rightarrow Cl_n$ be a $k$-hypermonogenic automorphic form on $\Gamma_{n-1}({\cal{I}})[N]$ and let ${\cal{R}}$ be a system of representatives in $\Gamma_{n-1}({\cal{I}})$ modulo $\Gamma_{n-1}({\cal{I}})[N]$.\\ Then for every $R \in {\cal{R}}$ the function $F(x;R):= \frac{\overline{(c_R x+d_R)}}{\|c_R x+d_R\|^{n+1-k}}f(R<x>)$ has a particular Fourier series representation of the form 
{\small{
\begin{eqnarray}
\label{khypfourier}
F(x;R) &=&  a(\underline{0};R) + \alpha(\underline{0};R) x_n^k \\
&+&\!\!\!\!\!  \sum_{\underline{\omega} \in \Lambda_{n-1}^{*}[N] \backslash \{\underline{0}\}} \!\!
x_n^{\frac{k+1}{2}} \Bigg[K_{\frac{k+1}{2}}(2 \pi \|\underline{\omega}\| x_n) - i e_n \frac{\underline{\omega}}{\| \underline{\omega}\|}  K_{\frac{k-1}{2}}(2 \pi \|\underline{\omega}\| x_n)\Bigg] \alpha(\underline{\omega};R) \nonumber \\
& & \quad\quad\quad\quad\quad \times 
 e^{2 \pi i <\underline{\omega},\underline{x}>}.\nonumber
\end{eqnarray}}}
Here, $a(0,R)$ and all $\alpha(\underline{\omega};R)$ are well-defined Clifford numbers from the Clifford sub algebra $Cl_{n-1}$. 
\end{theorem}
 
{\bf Remark.} In the monogenic case $(k=0)$, we have 
$$
K_{-\frac{1}{2}}(\|\underline{\omega}\| x_n) = K_{\frac{1}{2}}(\pi \|\underline{\omega}\| x_n) 
= \frac{\sqrt{\pi}}{\sqrt{2}\sqrt{\|\underline{\omega}\|x_n}} e^{-\|\underline{\omega}\|x_n}.$$ 
In this particular case, the factor  simplifies to 
$$
x_n^{\frac{1}{2}}\Bigg[K_{\frac{1}{2}}(\|\underline{\omega}\| x_n) - i e_n \frac{\underline{\omega}}{\| \underline{\omega}\|}  K_{-\frac{1}{2}}(\|\underline{\omega}\| x_n)\Bigg]
= \sqrt{\frac{\pi}{2}} \frac{1}{\sqrt{\|\underline{\omega}\|}}e^{- \pi \|\underline{\omega}\|x_n} (1-i e_n\frac{\underline{\omega}}{\| \underline{\omega}\|}).
$$
This is a scalar multiple of an idempotent in the Clifford algebra. In the monogenic case one obtains the well-known particular Fourier series representation, involving the monogenic plane wave exponential functions (see \cite{DSS,KraPhD}):
\begin{equation}
\label{fouriermonogenic}
F(x;R) = \alpha(\underline{0};R) + \sum\limits_{\underline{\omega} \in \Lambda_{n-1}^{*}({\cal{I}})[N] \backslash \{\underline{0}\}}  
(1 - i e_n \frac{\underline{\omega}}{\|\underline{\omega}\|}) \alpha(\underline{\omega};R) e^{2 \pi i <\underline{\omega},\underline{x}> - 2 \pi \|\underline{\omega}\| x_n}.
\end{equation}
In particular, one has the identity:
$$
\lim\limits_{x_n\rightarrow +\infty} F(x_n e_n;R) = \alpha(\underline{0};R),
$$
from which me may derive that monogenic cusp forms are exactly those monogenic automorphic forms whose first Fourier coefficient in the representation (\ref{fouriermonogenic}) vanishes. 

\par\medskip\par

In the more general $k$-holomorphic Cliffordian case we have the following

\begin{theorem}\label{hypercusp} (Fourier expansion of $k$-holomorphic Cliffordian cusp forms)\\ 
Suppose that $f: H^{+}(\mathbb{R} \oplus \mathbb{R}^n) \rightarrow Cl_n$ is a $k$-holomorphic Cliffordian automorphic form on $\Gamma_{n-1}({\cal{I}})[N]$. Suppose that for any system ${\cal{R}}$ of representatives of right cosets in $\Gamma_{n-1}({\cal{I}})$ modulo $\Gamma_{n-1}({\cal{I}})[N]$ each associated function $$F(x;R) := \frac{\overline{(cx+d)}}{\|cx+d\|^{n+1-k}} f(R<x>)$$ has a Fourier series representation of the particular form (\ref{khfourier}) where additionally all the Fourier images associated to the null frequency vanish identically. Then $f$ is a $k$-holomorphic Cliffordian cusp form of weight $(n-k)$ on $\Gamma_{n-1}({\cal{I}})[N]$.  
\end{theorem} 
{\bf Remark:} The Fourier images associated to the null frequency of $F(x;R)$, which explicitly reads $a(\underline{0};R) + \alpha(\underline{0};R)x_n^k + b(\underline{0};R) e_n x_n + \beta(\underline{0};R) e_n x_n^{k+1}$ vanish identically, if and only if the constants $a(\underline{0};R),\alpha(\underline{0};R), b(\underline{0};R)$ and $\beta(\underline{0};R)$ are equal to zero. 
\par\medskip\par

\subsection{The Petersson inner product and orthogonality}

In the context of the groups $\Gamma_{n-1}({\cal{I}})[N]$, the proper analogue of the Petersson inner product on the associated space of automorphic forms of weight $\alpha$ has the form (see also \cite{EGM90})
\begin{equation}
\label{petersson}
\langle f,g \rangle := \int\limits_{{\cal{F}}_{n-1}[N]} f(x) \overline{g(x)} x_n^{\alpha} d\nu(x)
\end{equation}
where again $d\nu(x) = x_n^{-n-1} dx_0 dx_1 \cdots dx_n$ stands for the invariant volume element of upper half-space associated to the hyperbolic metric. 
The integral is well-defined whenever at least one of the two functions $f$ and $g$ are cusp forms. It is independent of the choice of the fundamental domain, since $d\nu(x)$ is a $\Gamma_{n-1}({\cal{I}})[N]$-invariant measure and 
$$
f(M<x>) \overline{g(M<x>)} ({M<x>}_n)^{\alpha} = x_n^{\alpha} f(x) \overline{g(x)}
$$
holds for all $M \in \Gamma_{n-1}({\cal{I}})[N]$. 

In the monogenic case, the canonical attached weight factor is $\alpha=n$. In this particular case this inner product thus simplifies to
\begin{equation}
\label{peterssonmonogenic}
\langle f,g \rangle = \int\limits_{{\cal{F}}_{n-1}[N]} f(x) \overline{g(x)} x_n^{-1} dx_0 d x_1 \cdots dx_n.
\end{equation}
In the $k$-holomorphic Cliffordian case, this inner product has the form 
\begin{equation}
\label{peterssonkhypermonogenic}
\langle f,g \rangle = \int\limits_{{\cal{F}}_{n-1}[N]} f(x) \overline{g(x)} x_n^{-k-1} dx_0 d x_1 \cdots dx_n,
\end{equation}
from which we recover the monogenic case when putting $k=0$. 
   
\par\medskip\par

Theorem~\ref{kfnew} allows us now to establish: 

\begin{theorem}\label{ortho}
Suppose that $k < -1$ is an even integer and that $N \ge N_0({\cal{I}})$.  Let ${\cal{M}}_{k,n-1,N}$ be the space of $k$-holomorphic Cliffordian automorphic forms of weight $(n-k)$ with respect to $\Gamma_{n-1}({\cal{I}})[N]$. Let ${\cal{S}}_{k,n-1,N}$ be the subspace of $k$-holomorphic Cliffordian cusp forms of weight $(n-k)$ with respect to $\Gamma_{n-1}({\cal{I}})[N]$. 

Let further ${\cal{E}}_{k,n-1,N}$ be the space of the $k$-hypermonogenic Eisenstein series $\varepsilon_{k,n-1,N}(x)$ which we defined in (\ref{khypereisen}).  Then we have 
$$
{\cal{M}}_{k,n-1,N} = {\cal{E}}_{k,n-1,N} \oplus^{\perp} {\cal{S}}_{k,n-1,N}
$$
where orthogonality is understood in the sense of the Petersson inner product defined above. 
\end{theorem}
 
\begin{proof} Let ${\cal{V}}$ be an arbitrary but fixed chosen system of representatives of the right cosets of $\Gamma_{n-1}({\cal{I}})[N]$ modulo ${\cal{T}}_{n-1}({\cal{I}})[N]$. Furthermore, let ${\cal{R}}$ be an arbitrary system of representatives of right cosets in $\Gamma_{n-1}({\cal{I}})$ modulo $\Gamma_{n-1}({\cal{I}})[N]$. In the sequel we denote the closure of the fundamental domain $\Gamma_{n-1}({\cal{I}})[N]$ by $\overline{{\cal{F}}}_{n-1}$. We compute: 
{\small \begin{eqnarray*}
\label{orthorel}
\langle E, f \rangle &=& \int\limits_{\overline{{\cal{F}}}_{n-1}[N]} \sum\limits_{M \in {\cal{V}}}  \frac{\overline{(cx+d)}}{\|cx+d\|^{n+1-k}} \overline{f(x)} x_n^{n-k} d\nu(x)  
\end{eqnarray*}
\begin{eqnarray*}
& &= \sum\limits_{M \in  {\cal{V}}} \int\limits_{M<\overline{{\cal{F}}}_{n-1}[N]>} \frac{\overline{(cx+d)}}{\|cx+d\|^{n+1-k}}\overline{f(x)} x_n^{n-k} d\nu(x) \nonumber \\
& &= \int\limits_{x_n=0}^{x_n = +\infty} \sum\limits_{R \in {\cal{R}}}\Bigg[\int\limits_{
L_{n-1}[N]}  \overline{\frac{\overline{(c_Rx+d_R)}}{\|c_Rx+d_R\|^{n+1-k}} f(R<x>)} x_n^{-k-1} dx_0 dx_1 \cdots dx_{n-1}\Bigg] dx_n, \nonumber\\
& &= \int\limits_{x_n = 0}^{+\infty} \sum\limits_{R \in {\cal{R}}} \Big[ \int\limits_{L_{n-1}[N]} \overline{F(x;R)} x_n^{-k-1} dx_0 dx_1 \cdots dx_{n-1}\Big] dx_n,
\end{eqnarray*}}
where we put $F(x;R) := \frac{\overline{(c_R x+d_R)}}{\|c_Rx+d_R\|^{n+1-k}}f(R<x>)$. Since $f$ is $k$-holomorphic Cliffordian on upper half-space and $n$-fold periodic with respect to $\Lambda_{n-1}({\cal{I}})[N]$, it has a Fourier series expansion of the form (\ref{khypfourier}). In particular, for the first coefficient we have the identity   
\begin{eqnarray*}
& & \int\limits_{L_{n-1}} F(\underline{x}+x_n e_n;R) dx_0 dx_1 \cdots dx_{n-1}\\ &=& a(\underline{0};R) + \alpha(\underline{0};R) x_n^{k} + b(\underline{0};R) e_n + \beta(\underline{0};R) e_n x_n^{k+1}.
\end{eqnarray*}   
Hence,  
$$
\langle E,f\rangle = \int\limits_{x_n = 0}^{+\infty} \sum\limits_{R \in {\cal{R}}} \Big(a(\underline{0};R)+b(\underline{0};R)e_n x_n + \alpha(\underline{0};R)  x_n^{k}+\beta(\underline{0};R) e_n x_n^{k+1} \Big) x_n^{-k-1} dx_n.
$$
However, since $f$ is a $k$-holomorphic Cliffordian cusp form, we have $a(\underline{0};R)=0,b(\underline{0};R)=0, \alpha(\underline{0};R)=0$ and $\beta(\underline{0};R) = 0$ for all $R \in {\cal{R}}$. This proves the stated orthogonality result.  

\end{proof}
 
{\bf Remark}: This decomposition remains valid in the monogenic case $k=0$. As mentioned previously, in this case monogenic Eisenstein series can be introduced by (\ref{0hypereisen}), using the Hecke trick. However, we have not been able yet to prove the existence of non-vanishing monogenic cusp forms on $\Gamma_{n-1}({\cal{I}})[N]$.  Again, by applying Theorem~2 and Corollary~1 we can formulate a complete analogy of this decomposition statement for the cases where $k$ is an even positive integer in the context of $k+2$-holomorphic Cliffordian functions.

\end{document}